\input amstex
\input prepictex
\input pictex
\input postpictex
\magnification=\magstep1
\documentstyle{amsppt}
\TagsOnRight
\hsize=5.1in                                                  
\vsize=7.8in
\define\R{{\bold R}}
\define\cl{\operatorname {cl}}
\def\k{\bold k}

\vskip 2cm
\topmatter

\title Massey products and critical points \endtitle

\rightheadtext{}
\leftheadtext{}
\author  Michael Farber \endauthor
\address
School of Mathematical Sciences,
Tel-Aviv University,
Ramat-Aviv 69978, Israel
\endaddress
\email farber\@math.tau.ac.il
\endemail
\thanks{The research was supported by a grant from the 
Israel Academy of Sciences and Humanities and by
the Herman Minkowski Center for Geometry.} 
\endthanks
\abstract{In this paper we use cup-products and higher Massey products
to find topological lower bounds on the minimal number of geometrically distinct critical points
of any closed 1-form in a given cohomology class.}
\endabstract
\endtopmatter

\define\C{{\bold C}}

\define\Z{{\bold Z}}
\define\Q{{\bold Q}}

\redefine\c{{\frak c}}

\define\im{\operatorname{im}}



\define\E{\Cal E}

\def\<{\langle}
\def\>{\rangle}

\define\pd#1#2{\dfrac{\partial#1}{\partial#2}}

\def\part{\partial}

\NoBlackBoxes

\def\cat{\operatorname {cat}}

\def\fp{\Bbb F_p}

\documentstyle{amsppt}

\nopagenumbers

\heading{\bf \S 1. Introduction} \endheading

Let $X$ be a closed manifold and let $\xi\in H^1(X;\R)$ be a nonzero cohomology class. 
The well-known Novikov inequalities \cite{N1, N2} estimate the numbers of critical points of
different indices of any closed 1-form $\omega$ on $X$ lying in the class $\xi$, assuming that all the
singular points are non-degenerate in the sense of Morse. 
Novikov type inequalities were generalized in \cite{BF1} for closed 1-forms 
with more general singularities (non-degenerate in the sense of Bott). 
In \cite{BF2} an equivariant generalization of the Novikov inequalities was developed. 

Novikov inequalities have found important applications in symplectic topology, especially in the study of 
symplectic fixed points (Arnold's conjecture). 
Here we should mention the work of J.-C. Sikorav \cite{S}, Hofer - Salamon \cite{HS}, Van - Ono \cite{VO}, and most
recently the preprint of Eliashberg and Gromov \cite{EG}. 

In this paper we announce a new theorem, 
which provide topological restrictions on the number of geometrically
distinct critical points of closed 1-forms. We  impose no assumptions
on the nature of the critical points. Therefore, the result of this paper has the same 
relation to the classical Lusternik - Schnirelman - Frolov - Elsgoltz theory, as Novikov's theory has
to the classical Morse theory.

On any closed 
$n$-dimensional manifold, in any integral cohomology class $\xi$,
there always exists a closed 1-form with $\le n-1$ critical points (cf. Theorem 4.1).
Our Main Theorem 2.3 proves that on some manifolds
any closed 1-form has at least $n-1$ critical points. 

The main technical tool of our proof is {\it the deformation complex}, cf. \cite{F3}, \cite{F4}. 
It has advantages compared with the "Novikov's complex" (which uses formal power series). 

Different estimates on the number of critical points of closed 1-forms were recently suggested
in \cite{F2, F3}, where we used flat line bundles described by complex numbers, 
which are not Dirichlet units. 
There are examples, when the approach of this paper
gives stronger estimates than the approach of \cite{F2, F3}, although in some other
cases the situation is the opposite.

\heading{\bf \S 2. Main Theorem}\endheading

\subheading{2.1. Massey products} We will deal here with a special kind of higher Massey 
operations $d_r$ (where $r=1, 2, \dots$),
determined by a one-dimensional cohomology class $\xi\in H^1(X;\Z)$. The first operation $d_1: H^i(X;\k) \to 
H^{i+1}(X;\k) $ is the usual cup-product 
$$d_1(v)=(-1)^{i+1}v\cup\xi,\qquad\text{for}\quad v\in H^i(X;\k).\tag1$$
The higher Massey products $d_r: E_r^\ast \to E_r^{\ast+1}$ with $r\ge 2$, are defined as the differentials of 
a spectral sequence $(E_r^\ast, d_r)$, ($r\ge 1$)
with the initial term $E_1^\ast =H^\ast(X;\k)$ and the initial differential 
$d_1: E_1^\ast\to E_1^{\ast+1}$ given by (1). Each subsequent term $E_r^\ast$ is the cohomology of the previous
differential $d_{r-1}$: 
$$E_r^\ast=\ker(d_{r-1})/\im(d_{r-1}).\tag2$$
Traditionally, the following notation is used
$$d_r(v) = (-1)^{i+1}\cdot
\<v, \underset{\text{($r$ times)}}\to {\underbrace{\xi, \xi, \dots, \xi}}\>,\qquad
v\in H^i(X;\k).\tag3$$

This spectral sequence was first used in the context of the Novikov inequalities in paper
\cite{F}, page 46,
where we proved that the Novikov - Betti number $b_i(\xi)$ coincides with $\dim E_\infty^i$.
In \cite{N3} and in \cite{P} it is shown that the differentials $d_r$ are given by the Massey products
(3).

We describe this spectral sequence in full detail in \cite{F4}, cf. also \cite{F}, \cite{N3}, \cite{P},
\cite{F1}.

\subheading{Definition} A cohomology class $v\in H^i(X;\k)$ is said to be {\it $\xi$-surviver}
if $d_r(v)=0$ for all $r\ge 1$.

\subheading{2.2. Notation} 
Let $\k$ be a fixed algebraically closed field. The most important cases, which the reader should
keep in mind are $\k=\C$ or $\k$ being the algebraic closure of a finite field $\fp$.

We will consider flat $\k$-vector bundles $E$ over a compact polyhedron $X$. 
We will understand such bundles as
locally trivial sheaves of $\k$-vector spaces. The cohomology $H^q(X;E)$ will be understood as the
sheaf cohomology. 

A flat vector bundle is
determined by its monodromy -- linear representation of the fundamental group $\pi_1(X,x_0)$ 
on the fiber $E_0$ over the base point $x_0$, which is given by the parallel transport along loops.
For example, a flat $\k$-line bundle is determined by a homomorphism $H_1(X;\Z)\to \k^\ast$; here $\k^\ast$ is considered 
as a multiplicative abelian group.

The following Theorem is the main result of this paper.

\proclaim{2.3. Theorem} Let $X$ be a closed manifold and $\xi\in H^1(X;\Z)$ be an integral 
cohomology class. Suppose that there exists a nontrivial cup-product
$$v_1\cup v_2\cup \dots\cup v_m\ne 0,\tag4$$
where the first two classes $v_1\in H^{d_1}(X;\k)$ and $v_2\in H^{d_2}(X;\k)$ are $\xi$-survivors 
and for $i=3, \dots, m$ the classes $v_i\in H^{d_i}(X;E_i)$ belong to the cohomology of some flat 
$\k$-vector bundles
$E_i$ over $X$ with $d_i>0$. Then for
any closed 1-form $\omega$ on $X$ lying in class $\xi$,
$$\cat(S(\omega))\ge m-1\tag5$$
where $S(\omega)$ denotes the set of critical points of $\omega$ and $\cat$ denotes the Lusternik - Schnirelman category. In particular, 
the total number of geometrically distinct critical points satisfies
$$\#S(\omega)\ge m-1.\tag6$$
\endproclaim

Recall that a point $p\in X$ is a critical point of $\omega$ if $\omega_p =0$.

In case $\xi=0$ (when we study critical points of functions) the class $1\in H^0(X;\k)$ is a 
$\xi$-surviver. Hence in this case we may take $v_1=v_2=1$. This shows that in the case of functions
Theorem 2.3 is reduced to the usual Lusternik - Schnirelman inequality: 
$\cat(S(\omega))\ge \cl(X)+1,$
cf. \cite{DNF}.

\heading{\bf \S 3.  Examples}\endheading

\subheading{3.1. Detecting $\xi$-survivors} The following criterion allows us to show in some 
cases that a given cohomology class $v\in H^i(X;\k)$ is a $\xi$-surviver, where $\xi\in H^1(X;\Z)$.

{\it Suppose that we may realize $\xi$ by a smooth codimension one submanifold $V\subset X$,
having a trivial normal bundle, and we may realize the class $v$ by a simplicial cochain $c$,
such that the support of $c$ is disjoint from $V$. Then $v$ is a $\xi$-surviver.}

Indeed, this follows immediately from the definition of Massey operations.

This applies to the class $\xi$ itself and shows that it is a $\xi$-surviver, since we may realize it
by a cochain with support on a parallel copy of $V$. However,
this observation  is useless in producing nontrivial
products as in Theorem 2.3, because of the presence of another $\xi$-surviver, which kills the
whole product.

\subheading{3.2. Example} Let $X= \R\bold P^n\# (S^1\times S^{n-1})$ (the real 
projective space with a handle of index 1), and let $\xi\in H^1(X;\Z)$ be a class which 
restricts as the generator along the
handle and which is trivial on the projective space.
We may realize this class $\xi$ by a sphere $V=S^{n-1}$
cutting the handle. Let $\k=$ be the algebraic closure of $\Z_2$. We have a class
$v\in H^1(X;\k)$, which restricts as an obvious generator of $H^1(\R\bold P^n;\k)$ and which is trivial 
on the handle. Applying the above criterion, we see that $v$ is a $\xi$-surviver 
(since we may realize it by a chain with support on the projective space, 
i.e. disjoint from $V$). We obtain $v^n\ne 0$
and hence by Theorem 2.3 and closed 1-form $\omega$ in class $\xi$ has at least $n-1$ geometrically distinct critical points.

\subheading{3.3. Example} Let $\Sigma_g$ be a Riemann surface of genus $g>1$. Consider the 
classes $v_1, v_2, \xi\in H^1(\Sigma_g;\Z)$ which are Poincar\`e dual to the curves
shown in Figure 1.
\midinsert
$$\beginpicture
\setcoordinatesystem units <1.00000cm,1.00000cm>
\linethickness=1pt
\setlinear
%
%
\linethickness= 0.500pt
\circulararc 87.509 degrees from  3.302 21.145 center at  4.493 22.389
%
%
\linethickness= 0.500pt
\circulararc 50.506 degrees from  5.397 20.923 center at  4.443 18.987
%
%
\linethickness= 0.500pt
\circulararc 102.755 degrees from  7.588 21.336 center at  8.868 22.306
%
%
\linethickness= 0.500pt
\circulararc 46.809 degrees from  9.779 21.018 center at  8.858 18.891
%
%
\linethickness= 0.500pt
\circulararc 105.841 degrees from  4.381 22.796 center at  5.021 21.983
%
%
\linethickness= 0.500pt
\circulararc 102.072 degrees from  9.525 22.638 center at  9.782 21.638
%
%
\linethickness= 0.500pt
\ellipticalarc axes ratio  4.763:1.905  360 degrees 
	from 11.398 21.145 center at  6.636 21.145
%
%
\linethickness= 0.500pt
\ellipticalarc axes ratio  1.778:1.143  360 degrees 
	from  6.350 21.145 center at  4.572 21.145
\linethickness= 0.500pt
\setdots < 0.0953cm>
%
%
%
\plot	 4.413 22.860  4.572 22.701
 	 4.643 22.613
	 4.699 22.507
	 4.739 22.383
	 4.753 22.314
	 4.763 22.241
	 4.768 22.165
	 4.770 22.090
	 4.768 22.015
	 4.763 21.939
	 4.753 21.864
	 4.739 21.788
	 4.721 21.713
	 4.699 21.638
	 4.675 21.566
	 4.651 21.501
	 4.604 21.392
	 4.556 21.310
	 4.508 21.257
	 /
\plot  4.508 21.257  4.413 21.177 /
\linethickness= 0.500pt
%
%
%
\plot	 9.620 22.606  9.636 22.352
 	 9.638 22.288
	 9.635 22.221
	 9.628 22.153
	 9.616 22.082
	 9.600 22.010
	 9.580 21.935
	 9.555 21.859
	 9.525 21.780
	 9.492 21.704
	 9.455 21.632
	 9.414 21.565
	 9.370 21.503
	 9.272 21.394
	 9.160 21.304
	 /
\plot  9.160 21.304  8.922 21.145 /
%
%
\put{$v_2$} [lB] at  5.906 19.748
%
%
\put{$v_1$} [lB] at  4.477 23.209
%
%
\put{$\xi$} [lB] at  9.620 22.796
\linethickness=0pt
\putrectangle corners at  1.841 23.463 and 11.430 19.241
\endpicture
$$
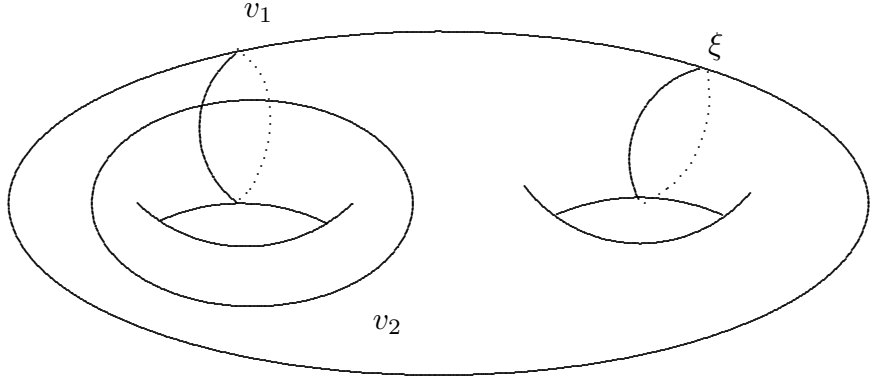
\centerline{Figure 1}
\endinsert
Then by 3.1, $v_1$ and $v_2$ are $\xi$-survivors. Since $v_1\cup v_2\ne 0$, we
obtain that any closed 1-form has at least 1 critical point. This however is a trivial corollary of the
Hopf's theorem. 

Let now $Y$ be an arbitrary closed manifold and let $X$ be $\Sigma_g\times Y$. 
Let $\xi'\in H^1(X;\Z)$ denote the class with $\xi'|_{\Sigma_g} =\xi$ and $\xi'|_Y=0$. 

We will show that:

{\it any closed 1-form on $X$ lying in class $\xi'$ has at least $\cl_{\k}(Y)$
critical points, where $\cl_{\k}(Y)$ is the cohomological cup-length of $H^\ast(Y;\k)$.}

Indeed, suppose that $u_j\in H^{d_j}(Y;\k)$ are some classes $j=1, 2, \dots, r= \cl_{\k}(Y)$ 
such that 
$d_j>0$ and 
$u_1\cup u_2\cup\dots \cup u_r\ne 0.$ Then we have $r+2$ cohomology classes 
$$v_1\times 1, \, v_2\times 1, \, 1\times u_1, \dots, 1\times u_r\in H^\ast(X;\k)$$
and their cup-product is nontrivial. Since the first two classes $v_1\times 1, v_2\times 1$
are $\xi'$-survivors (as easily follows from 3.1), we obtain our italicized
statement using Theorem 2.3.

This construction produces many examples with large numbers of critical points. 

For instance, one may take $Y$ being the real projective space  $\R\bold P^{n-2}$ and $\k$ being the algebraic closure of the field $\Z_2$. Then we obtain an $n$-dimensional manifold 
$$X= \Sigma_g\times \R\bold P^{n-2},\qquad g>1,$$
and the above arguments prove that {\it
any closed 1-form on $X$ has at least $n-1=\dim X - 1$ geometrically distinct critical points.}

\subheading{3.4. Formal spaces} If $X$ is a {\it formal space} \cite{DGMS} 
then all higher Massey products
vanish. We obtain in such a situation that 
{\it a class $v\in H^i(X;\Q)$ is a $\xi$-surviver if and only if $v\cup \xi=0$.} In this case the statement
of Theorem 2.3 simplifies. 

According to \cite{DGMS}, the vanishing of all higher Massey
products takes place in any compact complex manifold $X$
for which the $dd^c$-{\it Lemma} holds (for example, if $X$ is a K\"ahler or a Moi\v sezon space).
This vanishing of higher products directly follows from the diagram
$$\{\E^\ast_M, d\}\overset i\to\leftarrow \{\E^\c_M, d\}\overset \rho\to\to \{H_{d^c}(M),d\}$$
used in the first proof 
of the Main Theorem \cite{DGMS}, cf. page 270.

\subheading{3.5. Remark} {\it Theorem 2.3 becomes false if we allow products (4)
with only one $\xi$-surviver instead of two}. 

Indeed, let $X=S^1\times Y$ and let $\xi\in H^1(X;\Z)$
be the cohomology class of the projection $X\to S^1$. Then $\xi$ can be realized by a closed 
1-form without critical points (the projection). However,
products of the form $\xi\cup v_1\cup\dots \cup v_m$,
where $v_i$ are pullbacks of some classes of $Y$,
may be nontrivial for $m=\cl_\k(Y)$. Note that these cup-products have only one $\xi$-surviver.

\heading{\bf \S 4. Colliding critical points}\endheading

The following result shows that the examples described in \S 3 are best possible. 

\proclaim{4.1. Theorem} Let $X$ be a closed connected $n$-dimensional manifold, and let 
$\xi\in H^1(X;\Z)$ be a nonzero 
cohomology class. Then there exists a closed 1-form $\omega$
on $X$, realizing $\xi$, having at most $n-1$ critical points. \endproclaim

A proof can be obtained easily, using the method of Takens \cite{T}. It consists in colliding the Morse
critical points of the same dimension into one degenerate critical point.

Details and proofs will appear in \cite{F4}.

\enddocument